\documentclass[12pt]{article}
\usepackage{amsmath}
\usepackage{amssymb}
\usepackage{amsthm}

\font\elevenss=cmss11

\font\eightss=cmss8

\font\sixss=cmss8 at 6pt

\newfam\ssfam
\textfont\ssfam=\elevenss \scriptfont\ssfam=\eightss
 \scriptscriptfont\ssfam=\sixss

\theoremstyle{plain}
\newtheorem{thm}{Theorem}
\newtheorem{lem}[thm]{Lemma}

\theoremstyle{remark}

\def\ee{\epsilon}

\def\aa{{\bf a}}
\def\AA{{\bf A}}
\def\bb{{\bf b}}
\def\P{{\mathbb P}}

\def\Cox{\hfill \Box}
\def\C{{\mathbb C}}
\def\Z{{\mathbb Z}}

\makeatletter \def\romenumi{ \def\theenumi{\roman{enumi}}
\def\p@enumi{\theenumi} \def\labelenumi{(\@roman\c@enumi)}}
\makeatother

\date{Apr 20, 2009}

\title{Counting nondecreasing integer sequences that lie below a barrier}
\author{Robin Pemantle\\
University of Pennsylvania\\
Philadelphia, PA, 19104-6395\\
\texttt{pemantle@math.upenn.edu} \and
Herbert S. Wilf\\
University of Pennsylvania\\
Philadelphia, PA 19104-6395\\
\texttt{wilf@math.upenn.edu}}
\begin{document}
\maketitle

\begin{abstract}
Given a barrier $0 \leq b_0 \leq b_1 \leq \cdots$, let $f(n)$ be
the number of nondecreasing integer sequences $0 \leq a_0 \leq a_1
\leq \cdots \leq a_n$ for which $a_j \leq b_j$ for all $0 \leq j \leq n$.
Known formul\ae\ for $f(n)$ include an $n \times n$ determinant whose
entries are binomial coefficients (Kreweras, 1965) and, in the
special case of $b_j = rj+s$, a short explicit formula (Proctor, 1988, p.320).
A relatively easy bivariate recursion, decomposing all sequences
according to $n$ and $a_n$, leads to a bivariate generating function,
then a univariate generating function, then a linear recursion for
$\{ f(n) \}$.  Moreover, the coefficients
of the bivariate generating function have a probabilistic interpretation,
leading to an analytic inequality which is an identity for certain
values of its argument.
\end{abstract}

\section{Introduction}
We are given a sequence $\bb := b_0 \leq b_1 \leq \cdots$ of nonnegative
integers.  For another sequence $\aa := a_0 \leq a_1 \leq \dots$, we
write $\aa \preceq \bb$ if $a_i \leq b_i$ for all $i$.  Fix the upper
sequence $\bb$, which will be known as the \textit{barrier}.  Let
$f(n)$ denote the number of finite sequences $0 \leq a_0 \leq a_1 \leq
\cdots \leq a_n$ lying below the barrier.

The object of this paper is to find $f(n)$.  We first find
a generating function for $\{f(n)\}$, namely
\begin{equation}
\label{eq:krn1}
\sum_{n \geq 0} f(n)x^{n+1}(1-x)^{1+b_{n+1}} = 1 - (1-x)^{1 + b_0}  \, .
\end{equation}
The form of this function (mixed powers of $x$ and $1-x$) is
somewhat unusual but also it suggests, particularly if we replace
$x$ by $p$ and $1-x$ by $q$, that a probabilistic mechanism is operating
as well as the combinatorial one.
We then describe a family of random walks in the region below the
barrier, the family having the property that the $n$th term in
(\ref{eq:krn1}) above is the probability that a walk exits the region
under the barrier at exactly the $n$th step. However, it is possible
that with some positive probability a walk will never exit that
region, in which case the left-hand side of~\eqref{eq:krn1} will
be strictly less than the right-hand side.  This apparent paradox
is resolved by the observation that although (\ref{eq:krn1}) is
always true in the ring of formal power series (therefore determining
the numbers $f(n)$ uniquely by recursion), it is not necessarily true
as an analytic equality; in fact it will be true analytically
if and only if the probability is~1 that the walk eventually
collides with the barrier.

There is considerable overlap between this work and that of  \cite{gessel86}, in which Ira Gessel derived the functional equation (\ref{eq:krn1}) for the generating function by probabilistic, rather than combinatorial, means, and gave the random walk interpretation.

\section{Summary of results}
When we wish to emphasize the role
of $\bb$, we will write $f(\bb ; n)$ in place of $f(n)$.  Further we let
$\aa_{|m}$ denote the
truncation of $\aa$ to $(a_0 , \ldots , a_m)$, and extend the notation
to write $\aa_{|m} \preceq \bb$ when $a_i \leq b_i$ for $0 \leq i \leq m$.
Then formally, $f(\bb ; n)$ is the number of sequences $\aa_{|n}$ for
which $\aa_{|n} \preceq \bb$.  Equivalently, $f(\bb ; n)$ is the number of
integer partitions with $n$ parts whose Ferrers diagram is entirely
contained in the diagram for $\bb$, or the number of elements of
Young's lattice that lie below $\bb$.

An explicit formula for $f(\bb ; n)$ is known.  In fact, for the
more general problem of counting sequences of length $n$ lying
between two barriers, $\aa$ and $\bb$, Kreweras proved in
1965~\cite{kreweras65} that the number $f(\aa , \bb ; n)$ of
nondecreasing integer sequences between barriers $\aa$ and $\bb$
is given by the following $n \times n$ determinant:
\begin{equation} \label{eq:kreweras}
f(n) = \det{ \left(
   {b_i - a_j + n \choose i - j + n}_{i, j = 0,\dots, n} \right) } \, .
\end{equation}
Our first result is a recursion leading to a quadratic-time computation
for $f(n)$:
\begin{thm} \label{th:recur}
The numbers $\{ f(n) \}$ of nondecreasing integer sequences between
zero and $\bb$ satisfy the recurrence
\begin{equation} \label{eq:recur}
f(n) = (-1)^n {b_0 + 1 \choose n+1} + \sum_{0 \leq m < n}
   (-1)^{m+n+1}{1+b_{m+1} \choose n-m} f(m), \qquad(n=0,1,2,\dots) \, .
\end{equation}
\end{thm}

This theorem follows from the following formal power series identity.
\begin{thm} \label{th:FPS}
\begin{equation}
\label{eq:krn}
1 = (1-x)^{1 + b_0} + \sum_{n \geq 0} f(n)x^{n+1}(1-x)^{1+b_{n+1}}
\end{equation}
in the ring $\C [x]$.
\end{thm}

As a special case, we recover a formula in the case of a linear
barrier due to R.\ Proctor~\cite{proctor88}
(see also~\cite[(7.194)~in~exercise~7.101b]{EC2}): if $b_j = rj + s$ then
\begin{equation} \label{eq:proctor}
f(\bb ; n) = \frac{s+1}{n+1} {s + (n+1)(r+1) \choose n } \, .
\end{equation}

For convergence of the sum~\eqref{eq:krn} in the formal power series
ring, it is necessary, as is indeed that case, that there be only
finitely many summands for each monomal $x^k$.  One may ask whether
in fact~\eqref{eq:krn} holds as an analytic identity.  The probabilistic
interpretation of the summands in~\eqref{eq:krn}, leading to the
following result, is elaborated in the proof in Section~\ref{sec:probproofs}.
\begin{thm} \label{th:prob}
Fix a barrier $\bb$ and a number $p \in (0,1)$.  Let $q := 1-p$.
For each $n$ and each sequence $\aa := (a_0 \leq \cdots \leq a_n) \preceq \bb$
with length $\ell (\aa) = n+1$, define the weight
$$w(\aa) := \frac{p^{n+1} q^{b_{n+1} + 1}} {1 - q^{1+b_0}} \, .$$
Then
\begin{equation} \label{eq:ineq}
\sum_{\aa} w(\aa) \leq 1
\end{equation}
where the sum is over all finite sequences, that is, over all
$k > 0$ and all $\aa$ with $\ell (\aa) = k$.
Furthermore, if we let $\theta := \liminf b_n / n \in [0,\infty]$
then we can determine whether equality holds in~\eqref{eq:ineq}
in nearly all cases, as follows.
\begin{enumerate} \romenumi
\item if $\theta < q/p$ then equality holds;
\item if $\theta = q/p$ then equality holds provided
$|b_n - \theta n| = O(\sqrt{n})$;
\item if $\theta > q/p$ then equality fails.
\end{enumerate}
\end{thm}

\section{Combinatorial proofs}

First we will prove Theorem~\ref{th:FPS}, after which
Theorem~\ref{th:recur} follows almost immediately.

\noindent{\sc Proof of Theorem}~\ref{th:FPS}:  We begin by
identifying some elementary recursions.  Fix the barrier
sequence $\bb$ and let
$$c (m,j) := f(\bb_{|m} ; j)$$
count the class $\Gamma (m,j)$ of sequences $\aa \preceq \bb$ of
length $m+1$ for which $a_m = j$.  We define $c(m,j) := 0$ if $m < 0$
or if $j > b_m$.  Partitioning $\Gamma (m,j)$ according to the
value of $a_{m-1}$ gives a disjoint union $\bigcup_i \Gamma (m-1 , i)$.
The set $\Gamma (m-1,i)$ is empty when $i > j$, whence
\begin{equation} \label{eq:gamma}
c(m,j) = \sum_{i=0}^j c(m-1 , i) \, .
\end{equation}
In particular, setting $j = b_m$,
\begin{equation} \label{eq:f(n)}
f(m-1) = c(m , b_m) \, .
\end{equation}
Next, we partition $\Gamma (m,j)$ into $A \cup B$ where $A$ is the
set of sequences $\aa$ with $a_{m-1} = a_m = j$ and $B$ is the
set of sequences with $a_{m-1} < a_m = j$.  The set $B$ is in
bijection with the set $\Gamma (m,j-1)$ via the map that changes
$a_m$ from $j$ to $j-1$ and fixes $a_i$ for $i < m$.  Clearly,
$|A| = c(m-1,j)$.  This implies the relation
\begin{equation} \label{eq:c}
c(m , j) = c(m-1 , j) + c(m , j-1)
\end{equation}
for every $m > 0$ and $0 < j \leq b_m$.  Checking the boundary
case $j = 0$, we have $c(m,0) = 1$ for $m \geq 0$ and zero
otherwise, so~\eqref{eq:c} holds for $j=0$ as long as $m \neq 0$.
When $m=0$, we have $c(0,j) = 1$ if $0 \leq j \leq b_0$ and
zero otherwise, whence~\eqref{eq:c} holds for $m=0$ as long
as $j \notin \{ 0 , b_0 + 1 \}$.  When $j \geq b_m + 2$, all
the terms of~\eqref{eq:c} are zero and the relation holds vacuously.
Finally, for any $m \geq 0$ and $j = b_m + 1$ we have
\begin{equation} \label{eq:not}
c(m , j) - c(m-1 , j) - c(m , j-1) = - c(m , b_m) = - f(m-1)
\end{equation}
by~\eqref{eq:f(n)}.  Define a bivariate generating function
$$C(x,t) := \sum_{m,j \geq 0} c(m,j) x^m t^j \, .$$
The relations~\eqref{eq:c} and exceptions~\eqref{eq:not} imply that
\begin{equation} \label{eq:kernel}
(1-x-t) C(x,t) = 1 - t^{1+b_0} -
   \sum_{m > 0} f(m-1) x^m t^{1+b_m} \, .
\end{equation}
The kernel method (see, e.g.,~\cite{B-MP,flajolet-sedgewick-anacomb})
suggests the substitution $t = 1 - x$.  On both sides of~\eqref{eq:kernel}
the power of $t$ appearing
in any monomial $x^m t^j$ is bounded by $b_m +1$, hence the substitution
is valid in the ring of formal power series and yields
$$0 = 1 - (1-x)^{1 + b_0} - \sum_{m > 0} f(m-1) x^m (1-x)^{1 + b_m} \, .$$
With $m = n+1$, this is Theorem~\ref{th:FPS}.
$\Cox$

\noindent{\sc Proof of Theorem}~\ref{th:recur}: For $k \geq 0$, the
coefficient of $x^{k+1}$ on the right-hand side of~\eqref{eq:krn}
is known to vanish.  But this coefficient is equal to
$$(-1)^{k+1} {b_0 + 1 \choose k+1} + \sum_{m=0}^k f(m) (-1)^{k-m}
   {b_{m+1} + 1 \choose k-m} \, .$$
Solving for $f(k)$ gives
$$f(k) = - \left [ (-1)^{k+1} {b_0 + 1 \choose k+1} + \sum_{m=0}^{k-1}
   f(m) (-1)^{k-m}  {b_{m+1} + 1 \choose k-m} \right ]$$
which is Theorem~\ref{th:recur} with the variable $k$ in place of $n$.
$\Cox$

To prove~\eqref{eq:proctor}, let $F(x) := \sum_{n \geq 0} f(n) x^n$
be the generating function for $\{ f(n) \}$.  Again, substitute
$t = 1-x$ in~\eqref{eq:kernel}; the left-hand side is again zero,
while the choice of $b_n = rn+s$ makes the right-hand side into
\begin{equation} \label{eq:F}
1 - (1-x)^{s+1} - x (1-x)^{r+s+1} F(x(1-x)^r) \, .
\end{equation}
There is a unique formal power series $X(y)$ with no constant term
such that $X(y) (1-X(y))^r = y$.  From~\eqref{eq:F} we get
$$F(x(1-x)^r) = \frac{1 - (1-x)^{s+1}}{x(1-x)^{r+s+1}} \, .$$
Composing formally with $X$ we obtain
$$F(y) = \frac{1 - (1 - X(y))^{s+1}}{y (1-X(y))^{s+1}} \, .$$
Thus,
\begin{equation} \label{eq:lagrange}
f(n) = [y^n] F(y) = [y^{n+1}] (y F(y)) = [y^{n+1}] \left (
   \frac{1}{(1 - X(y))^{s+1}} - 1 \right ) \,
\end{equation}
which we may now evaluate via Lagrange inversion.  The following form
of the Lagrange inversion formula may be found in~\cite{wilf-GFology}.
\begin{lem} \label{lem:lagrange}
Let $\phi$ be formal power series in $x$ with $\phi (0) = 1$.
Then there is a unique formal power series $x = x(y)$ satisfying
$x = y \phi (x)$.  Further, if this series $x(y)$ is substituted
into another formal power series $H$, then the resulting series
satisfies
$$[y^n] H(x(y)) = \frac{1}{n} [x^{n-1}] \{ H' (x) \phi (x)^n \} \, .$$
$\Cox$
\end{lem}
When $\displaystyle{\phi (x) = \frac{1}{(1-x)^r}}$ then the series
$x(y)$ is the series $X(y)$ above.  In this case,
$$H(x(y)) = \frac{1}{(1 - X(y))^{s+1}} - 1$$
is the function on the right-hand side of~\eqref{eq:lagrange}.
Lagrange inversion with $n+1$ in place of $n$ gives
$H' (x) = (s+1)/(1-x)^{s+2}$, whence
\begin{eqnarray*}
f(n) & = & [y^{n+1}] \left ( \frac{1}{(1 - X(y))^{s+1}} - 1 \right ) \\
& = & \frac{1}{n+1} [x^n] \left \{ \frac{s+1}{(1-x)^{s+2+(n+1)r}} \right \} \\
& = & \frac{s+1}{n+1} {s + (n+1)(r+1) \choose n} \, ,
\end{eqnarray*}
proving~\eqref{eq:proctor}.

\section{Probabilistic proofs} \label{sec:probproofs}

Let $L$ be the set of points in $\Z^2$ defined by $L :=
\{ (i-1,j) : 0 \leq i ,  0 \leq j \leq b_i \}$.  Fix $0 < p < 1$
and let $\Omega$ be the space of infinite sequences of 0's and 1's,
equipped with the product measure $\P$ making each coordinate
0 with probability~$p$ and 1 with probability $q := 1-p$.
With each $\omega \in \Omega$ we associate a lattice path
beginning at the location $(-1,0)$, moving upward on step
$k$ when $\omega(k) = 1$, and moving right on step $k$ when
$\omega (k) = 0$, for each $k=0,1,2,\dots$.  If we let $S(k) :=
S(k , \omega) = \sum_{j=0}^k \omega (j)$, then a formal definition
of the induced path is the sequence $\{ (X(k),S(k)) : k \geq 0 \}$
of random vectors where $X(k) := k - S(k)$.
If $\omega$ begins with a block of more than $b_0$ 1's then the walk
will be outside of $L$ when it takes its first step to the right, so
to such an $\omega$ we associate the empty path.

Let $\tau (\omega)$ be the stopping time defined by $\tau :=
\inf \{ k \geq 0 : S(k) > b_{1+X(k)} \}$.  In other words, it is the first
time $k$ that $X(k) \notin L$ (the barrier is exceeded).  For each path,
we now pass to the subsequence corresponding to the locations after
moves to the right.  Formally, define the random variable $M$ by
$$M := \sup \{ X(k) : k \leq \tau \}$$
to be the farthest right extent of the path before exiting the
barrier and define a random sequence $\AA$ of length $M+1$ by
$$A_i := \min \{ j : X(k) = (i , j) \mbox{ for some } k \leq \tau \} \, .$$
It is possible that $M$ is infinite (the barrier is never reached).
It is also possible that $M = -1$, if initially there are $b_0+1$
upward moves, in which case, as remarked earlier, we set $\AA := \emptyset$, the empty path.

\noindent{\sc Proof of Theorem}~\ref{th:prob}:
To prove the inequality~\eqref{eq:ineq}, it suffices to observe that
for every sequence $\aa$ of length $m+1$, the probability that $M = m$ and
$\AA = \aa$ is equal to $p^{m+1} q^{b_{m+1} + 1}$.  Indeed, the event
$\{ \AA = \aa \}$ requires a specific sequence of values of
$\omega (k)$ for $1 \leq k \leq m + b_m + 2$, namely, upward
moves to height $a_0$, then a move to the right, then upward moves
to height $a_1$, then moves to the right, etc., ending with a move
to the right that ends at $(m , a_m)$, followed by upward moves
to height $b_{m+1} + 1$; the total number of rightward moves is $m+1$
(remember, we started at $(-1,0)$) and the total number of upward
moves is $b_{m+1}+1$, verifying the formula
$$\P (\AA = \aa) = p^{m+1} q^{b_{m+1} + 1} \, .$$
The event that $\AA = \emptyset$ has probability $q^{b_0 + 1}$.
Conditioning on this not occurring, we have
$$1 \geq \P (M < \infty | M > -1) = \frac{
   \sum_{\aa} p^{\ell (a)} q^{b_{\ell (a)}+1}}
   {1 - q^{b_0 + 1}}
   = \sum_{\aa} w(\aa) \, ,$$
proving~\eqref{eq:ineq}.

Evidently, equality holds if and only if $\P (M = \infty) = 0$.
This is the well known problem of whether a random walk can remain
forever on one side of a barrier.  An exact summability criterion
is known for this under some regularity assumptions on the barrier.
For example, in~\cite[(0.13)~and~Theorems~1 and~3]{erdos42}, Erd\"os
proves a summability criterion in the case where $p=1/2$ and
$n^{-1/2} (b_n - n/2)$ is nondecreasing.  The earliest version of such
a test in the continuous time case is due to Petrowsky~\cite{petrowsky35}.

Our theorem does not require results as sharp as these.  It suffices
to observe that if $(Y_n , Z_n)$ are the coordinates of $X_n$ then
the strong law of large numbers implies that $Z_n / Y_n \to q/p$
almost surely.  This implies that
$\P (Z_n \geq b_{Y_n} \mbox{ infinitely often}) = 0$
when $\theta > q/p$, which implies that $\P (Z_n \geq b_{Y_n}
\mbox{ for some } \, n ) < 1$, which implies $\P (M = \infty) > 0$.
Conversely, if $\theta < q/p$, then the strong law of large numbers implies
$\P (Z_n > (\theta + \ee) Y_n \mbox{ for sufficiently large } \, n) = 1$,
which implies $\P (M = \infty) = 0$.  Finally, if $\theta = q/p$,
the law of the iterated logarithm~\cite[Theorem~(9.7)]{durrett}
gives $Z_n \geq \theta Y_n + C \sqrt{Y_n}$ infinitely often
almost surely for any $C$, which implies $\P (M = \infty) = 0$
under the assumption $|b_n - \theta n| \leq C n^{1/2}$.  This completes
the proof of Theorem~\ref{th:prob}.
$\Cox$

\section{A question, and some acknowledgments}
We have not been able to generalize this combinatorial/probabilistic
method to the situation where there is a lower, as well as an upper,
barrier.  Nonetheless the result of Kreweras cited above suggests
that this may be possible.

Our thanks go to Mireille Bousquet-M\'elou and Richard Stanley for citations
to earlier work on this problem, and to Davar Khoshnevisan for citations to
summability criteria for random walks.

\newpage

\bibliographystyle{alpha}

\end{document}